\documentclass[12pt,a4paper]{article} 
\usepackage{amsmath}
\usepackage{amsfonts}
\usepackage{amssymb}
\usepackage{cite}
\usepackage{graphics}
\usepackage{amssymb,amscd}
\usepackage[utf8]{inputenc}
\usepackage[T2A]{fontenc}
\usepackage{tikz-cd}
\usepackage[russian,english]{babel}

\usepackage{amscd}
\usepackage{yfonts}
\usepackage[normalem]{ulem}
\usepackage{esint}
\usepackage{amsmath}
\usepackage{amsfonts}
\usepackage{amssymb}
\usepackage{hyperref}

\usepackage{graphics}
\usepackage{amssymb,amsmath,mathrsfs,amsthm,amscd}

\theoremstyle{plain}

\newtheorem{theo}{Theorem}[section]
\newtheorem{prop}[theo]{Proposition}

\theoremstyle{definition}

\theoremstyle{remark}

\numberwithin{equation}{section}

\DeclareMathOperator{\tr}{tr}

\DeclareMathOperator{\Scal}{Scal}
\DeclareMathOperator{\Ric}{Ric}
\DeclareMathOperator{\R}{R}
\DeclareMathOperator{\Vol}{vol}
\DeclareMathOperator{\dvol}{dvol}

\title{Heat trace asymptotics for wedge-like singularity.}
\author{Asilya Suleymanova}
\date{October 2017}

\begin{document}
\maketitle

\begin{abstract}
In this note we consider a heat trace expansion on a manifold with wedge-like singularity. We show that there are two terms in the expansion that contain information about the presence of the singularity, namely the logarithmic term $ct^{-1/2}\log t$ and the half power term $+b t^{-1/2}$. We also give a geometric expression for $c$.
\end{abstract}

\tableofcontents

\begin{section}{Geometric setup}
Consider a Riemannian manifold, $(M,g)$, of dimension $m$ such that $M=M_1\cup U$, where $M_1$ is a compact manifold with boundary and $U=(0,1)\times N$ has a fibration $(0,1)\times F\to U\to(0,1)\times \Sigma$, where $N,\Sigma$ are compact. This fibration induces a vertical bundle of tangents to the fibres $T_V N\subset TN$ and a horizontal cotangent subbundle $T^*_HN$ of cotangents annihilating $T_VN$. We take a complementary tangent subbundle $T_HN\subset TN$, inducing a vertical $T^*_VN$ that annihilates $T_HN$ such that $TN=T_VN\oplus T_HN$.

The metric on $U$ is
$$
g_U=dr^2+g_H+r^2g_V(s),\,\,\,\,\,\,\, r\in(0,1),
$$
where $g_H$ is induced from a metric on $\Sigma$ and $g_V(s)$ is a family of metrics on the fibres, $s\in\Sigma$.

Assume now that $\dim\Sigma=1$, i.e. it is a circle. Assume furthermore that metric on the fibres does not depend on $s$, i.e. $g_V(s)=g_V$.
$$
g_{\text{wedge}}=dr^2+d\theta^2+r^2g_V,\,\,\,\,\,\,\, r\in(0,1),
$$
We call such singularity a wedge-like singularity.
\end{section}

\begin{section}{Curvature tensors near the singularity}
In this section we give explicit formulas for the curvature tensors in $(U,g_{\text{wedge}})$ in terms of the curvature tensors on fibers. For the latter we use the classical tensor notations and for the tensors on $(U,g_{\text{wedge}})$ we use classical notations but with tilde.

Let $(x^1,\dots,x^{m-2})$ be local coordinates on a fiber $F$ and \\$p=(r,\theta,x^1,\dots,x^{m-2})\in U$. For $i,j\in\{r,\theta,1,\dots,m-2\}$ denote by $\tilde{g}_{ij}$ the components of the metric tensor on $g_{\text{wedge}}$, and by $g_{ij}$ for $i,j\in\{1,\dots,m-2\}$ the components of the metric tensor $g_V$. Then
$$
\tilde{g}_{rr}=\tilde{g}_{\theta\theta}=1,
$$
$$
\tilde{g}_{ir}=\tilde{g}_{ri}=\tilde{g}_{j\theta}=\tilde{g}_{\theta j}=0,\,\,\text{ for }i\in\{\theta,1,\dots,m-2\} \text{ and }j\in\{r,1,\dots,m-2\}
$$
and
$$
\tilde{g}_{ij}=r^2g_{ij}, \,\,\text{ for }i,j\in\{1,\dots,m-2\}.
$$

Then
\begin{align*}
\tilde{g}_{ij,k}(r,\theta,x)=
\begin{cases}
2rg_{ij}(x),\,\,\text{ if }k=r,\\
0 \,\,\,\,\,\,\,\,\,\,\;\;\;\;\;\;\;\text{ if }k=\theta,\\
r^2g_{ij,k}(x)\,\,\text{ otherwise.}
\end{cases}
\end{align*}

The Christoffel symbols are of course
$$
\tilde{\Gamma}^i_{jk}=\frac12 \tilde{g}^{im}\left(\tilde{g}_{mj,k}+\tilde{g}_{mk,j}-\tilde{g}_{jk,m}\right),
$$
and now we express them in terms of the Christoffel symbols $\Gamma^i_{jk}$ and the metric tensor $g_{ij}$.

Let $i=r$, then

$$
\tilde{\Gamma}^r_{jk}=
\begin{cases}
0,\,\,\text{ if }j=r \text{ or }k=r \text{ or }j=\theta \text{ or }k=\theta,\\
-rg_{jk}\,\,\text{ otherwise.}
\end{cases}
$$

Let $i\notin\{r,\theta\}$ and $j=r$, then

$$
\tilde{\Gamma}^i_{rk}=
\tilde{\Gamma}^i_{kr}=
\begin{cases}
0,\,\,\text{ if }k=r \text{ or }i=\theta \text{ or }k=\theta,\\
r^{-1}\delta^i_k\,\,\text{ otherwise.}
\end{cases}
$$

Assume that $i,j,k\notin\{r,\theta\}$, then $\tilde{\Gamma}^i_{jk}=\Gamma^i_{jk}$. Finally $\tilde{\Gamma}^\theta_{jk}=\tilde{\Gamma}^i_{\theta k}=\tilde{\Gamma}^i_{k\theta}=0$.

The scalar curvature $\tilde{\Scal}(p)$ can be expressed in the following way

\begin{equation}\label{scalar curvature}
\begin{split}
\tilde{\Scal}(p)
=&\tilde{g}^{ij}\left(\tilde{\Gamma}^m_{ij,m}-\tilde{\Gamma}^m_{im,j}
+\tilde{\Gamma}^l_{ij}\tilde{\Gamma}^m_{ml}-\tilde{\Gamma}^l_{im}\tilde{\Gamma}^m_{jl}\right)\\
=&r^{-2}\big(\Scal(x)-(m-2)(m-3)\big),
\end{split}
\end{equation}
where $p=(r,\theta,x)\in U$, $x\in F$ and $\tilde{\Scal}(p)$ is the scalar curvature on $(U,g_{\text{wedge}})$, and $\Scal(x)$ is the scalar curvature on $(F,g_V)$.

If $i,j,k,l\notin\{r,\theta\}$, we obtain

\begin{equation}\label{Riemann curvature}
\tilde{\R}_{ijkl}(p)=r^{-2}\big(\R_{ijkl}(x)-g_{ip}(x)g_{jm}(x)(\delta^p_k\delta^m_l-\delta^p_l\delta^m_k)\big).
\end{equation}

Similarly, for the tensor Ricci
\begin{equation}\label{Ricci curvature}
\tilde{\Ric}_{ij}(p)=r^{-2}\left(\Ric_{ij}(x)-(m-3)g_{ij}(x)\right).
\end{equation}

\end{section}

\begin{section}{Resolvent expansion}

Let $\Delta$ be the Friedrichs extension of the Laplace-Beltrami operator on $U$
$$
\Delta=-\partial_r^2+\Delta_\Sigma+r^{-2}A(s),
$$
where $\Delta_\Sigma$ is the Laplace operator on the singular stratum $\Sigma$, $A(s), s\in\Sigma$ is a smooth family of elliptic operators on functions on fibers $F$. In the case of the wedge-like singularity, we have $\Delta_\Sigma=-\partial_\theta^2$, and $A(s)=\Delta_F+\left(\frac{m-2}{2}\left(\frac{m-2}{2}-1\right)\right)$, where $\Delta_F$ is the Laplace-Beltrami operator on the fiber $F$.

Let $\varphi(r,\theta)$ be a cutoff function supported near $r=0$, by \cite[Theorem~5.2]{BS}, for $d>m/2$ the resolvent $(\Delta+z^2)^{-d}$ is trace class and the following is true
$$
\tr(\varphi(r,\theta)(\Delta+z^2)^{-d}))=\int_0^{\infty}\varphi(r,\theta)\sigma(r,rz)dr,
$$ 
where $\sigma(r,rz):=\sigma(r,\theta,rz):=\tr_{L^2(F)}(\Delta+z^2)^{-d}.$

From the local heat kernel expansion of heat kernel away from the singularity we know, see \cite[(3.6), (3.17)]{S}, that
\begin{align*}
&\tr_{L^2(F)}(\Delta+\zeta^2/r^2)^{-d}
\sim_{\zeta\to\infty}\\
&(4\pi)^{-\frac{m}{2}}\sum_{j=0}^{\infty}
(\zeta/r)^{-2d+m-2j}\frac{\Gamma(-\frac{m}{2}+d+j)}{(d-1)!}\int_Fr^{m-2}u_{j}(p)\dvol_F.
\end{align*}
Above $u_{j}(p)$ are polynomials in the curvature tensors at $p\in M$ that come from the local heat kernel expansion away from the singularity, in particular $u_0(p)\equiv1$, $u_1(p)=\Scal(p)$ and 
$$
u_2(p)=\frac{1}{360}\left(12\Delta\Scal(p)+5\Scal(p)^2-2\lvert\Ric(p)\rvert^2+2\lvert\R(p)\rvert^2\right).
$$

Hence
\begin{equation}\label{sigma_expansion}
\sigma(r,\zeta)
\sim_{\zeta\to\infty}\sum_{j=0}^{\infty}\zeta^{-2d+m-2j}\sigma_{j}(r),
\end{equation}
where

\begin{equation}\label{sigma expansion terms}
\begin{split}
\sigma_{j}(r)
=&(4\pi)^{-\frac{m}{2}}r^{2d-m+2j}\int_Fr^{m-2}u_{j}(p)\dvol_N\frac{\Gamma(-\frac{m}{2}+d+j)}{(d-1)!}\\
=&(4\pi)^{-\frac{m}{2}}r^{2d-2+2j}\int_Fu_{j}(p)\dvol_N\frac{\Gamma(-\frac{m}{2}+d+j)}{(d-1)!}.
\end{split}
\end{equation}

We can apply the Singular Asymptotics Lemma to obtain an asymptotic expansion of the resolvent trace, see \cite[p. 287]{BS},

\begin{prop}\label{applySAL}
\begin{align}
\int_0^{\infty}\varphi(r,\theta)\sigma(r,\zeta)dr
\sim&\sum_{l=0}^{\infty}z^{-l-1}\frac{1}{l!}\fint_0^{\infty}\int_{\Sigma}\zeta^l\partial^l_r\bigg(\sigma(r,\zeta)\varphi(r,\theta)\bigg)|_{r=0}d\theta d\zeta\\
+&\sum_{j=0}^{\infty}\fint_0^{\infty}\int_{\Sigma}\sigma_{j}(r)(rz)^{-2d+m-2j}\varphi(r)d\theta dr\\
+&\sum_{l=\frac{m}{2}-d+1}^{\infty}z^{-2d+m-2l}\log z\int_{\Sigma}\frac{\partial^{2d-m+2l-1}_r\bigg(\sigma_{l}(r)\varphi(r,\theta)\bigg)|_{r=0}}{(2d-m+2l-1)!}d\theta.
\end{align}
\end{prop}

Assume $\varphi(r,\theta)\equiv1$ near $r=0$ and consider the last sum in the expansion

\begin{align*}
L:=&\sum_{l=\frac{m}{2}-d+1}^{\infty}z^{-2d+m-2l}\log z\int_{\Sigma}\frac{\partial^{2d-m+2l-1}_r\bigg(\sigma_{l}(r)\varphi(r,\theta)\bigg)|_{r=0}}{(2d-m+2l-1)!}d\theta\\
=&\sum_{j=\frac{m}{2}-d+1}^{\infty}
z^{-2d+m-2j}\log z\frac{\Gamma(-\frac{m}{2}+d+j)\Vol(\Sigma)}{(d-1)!(2d-m+2j-1)!}\times\\
&\times(4\pi)^{-\frac{m}{2}}
\partial^{2d-m+2j-1}_r\bigg(r^{2d-2+2j}\int_F u_j(p)\dvol_N\bigg)|_{r=0}.
\end{align*}

To continue the computation we note that for $j\geq m/2$ we have

\begin{equation*}
\begin{split}\partial^{2d-m+2j-1}_r\bigg(r^{2d-2+2j}\int_F u_j(p)\dvol_N\bigg)|_{r=0}\\
=\partial^{2j-m+1}_r\bigg(r^{2j}\int_F u_j(p)\dvol_N\bigg)|_{r=0},
\end{split}
\end{equation*}
this makes sense only for $j\geq\frac{m-1}{2}$,
therefore

\begin{equation*}
\begin{split}
L=&\sum_{j\geq\frac{m-1}{2}}
z^{-2d+m-2j}\log z\frac{\Gamma(-\frac{m}{2}+d+j)\Vol(\Sigma)}{(d-1)!(2d-m+2j-1)!}\times\\
&\times(4\pi)^{-\frac{m}{2}}
\partial^{2j-m+1}_r\bigg(r^{2j}\int_F u_j(p)\dvol_F\bigg)|_{r=0}\\
=&\sum_{l=0}^{\infty}
z^{-2d-2l+1}\log z\frac{\Gamma(d+l-\frac12)\Vol(\Sigma)}{(d-1)!(2d+2l-2)!}\times\\
&\times(4\pi)^{-\frac{m}{2}}
\partial^{2l}_r\bigg(r^{m+2l-1}\int_F u_{\frac{m-1}{2}+l}(p)\dvol_F\bigg)|_{r=0}\\
=&z^{-2d+1}\log z\frac{\Gamma(d-\frac12)\Vol(\Sigma)}{(d-1)!(2d-2)!}\times\\
&\times(4\pi)^{-\frac{m}{2}}
\bigg(r^{m-1}\int_F u_{\frac{m-1}{2}}(p)\dvol_F\bigg)|_{r=0}.
\end{split}
\end{equation*}
The last equality is due to the fact that $r^{m-1+2l}u_{\frac{m-1}{2}+l}(p)$ is a smooth function with respect to $r$.

By \cite[p.47]{S}, the logarithmic part in the heat trace expansion coming from this term is

\begin{align}\label{logarithmic term}
-(4\pi)^{-\frac{m}{2}}
t^{-\frac12}\log t\times
\frac{1}{2}\Vol(\Sigma)\bigg(r^{m-1}\int_F u_{\frac{m-1}{2}}(p)\dvol_F\bigg)|_{r=0}.
\end{align}

\end{section}

\begin{section}{Heat trace expansion}

From Proposition~\ref{applySAL} and (\ref{logarithmic term}) we obtain
\begin{align}
\tr e^{-t\Delta}\sim_{t\to0+}
(4\pi t)^{-\frac{m}{2}}
\sum_{j=0}^{\infty}\tilde{a}_jt^j
+b t^{-1/2}
+ct^{-1/2}\log t,
\end{align}
where
\begin{align*}
\tilde{a}_j=
\begin{cases}
\int_Mu_j\dvol_M \text{ for } j\leq m/2-1, \\
\fint_Mu_j\dvol_M \text{ for } j>m/2-1.
\end{cases}
\end{align*}
and
$$
c=\frac{1}{2}\Vol(\Sigma)\bigg(r^{m-1}\int_F u_{\frac{m-1}{2}}(p)\dvol_F\bigg)|_{r=0},
$$
$b$ is to be computed.

Note that if $m$ is even, the logarithmic term is zero.

\begin{subsection}{Logarithmic term in low dimensional cases}

Let $\dim M=3$, then $c=\frac{1}{12}\Vol(\Sigma)\Scal(x)=0$, so the logarithmic term is always zero in this case.

Let $\dim M=5$, then 
\begin{align*}
c
=&\frac{1}{2}\Vol(\Sigma)\int_F\left(\frac{1}{360}\left(5\Scal(p)^2-2\lvert\Ric(p)\rvert^2+2\lvert\R(p)\rvert^2\right)\right)\dvol_F\\
=&-\frac{\Vol(\Sigma)}{720}(4\pi)^{-2}\int_F\Big(3(\Scal(x)-6)^2+6(\Ric_{ij}(x)-2g_{ij})^2\Big)\dvol_F.
\end{align*}
The above expression is equal to zero if and only if $\Scal(x)\equiv6$ and $\Ric_{ij}(x)=2g_{ij}(x)$; equivalently if and only if the sectional curvature of the fiber $(F,g_V)$ is equal to one. Therefore the logarithmic term in this case is equal to zero if and only if the fiber $(F,g_V)$ is isometric to a spherical space form.

\end{subsection}

\end{section}


\begin{thebibliography} {3}
\bibitem[BS]{BS} J.\,Br\"uning, R.\,Seeley, {\em The expansion of the resolvent near a singular stratum of conical type}, J.~Funct.~Anal. 95 (1991), 255--290.
\bibitem[S]{S} A.\,Suleymanova, {On the spectral geometry of manifolds with conic singularities}, PhD thesis, \url{https://edoc.hu-berlin.de/handle/18452/19097} (2017).
\end{thebibliography}
\end{document}